\documentclass[11pt]{article}

\usepackage{graphicx}
\usepackage{amsmath}
\usepackage{amsfonts}
\usepackage{fancyhdr} 
\usepackage{mathrsfs}

\usepackage{wrapfig}

\usepackage{bbold}

\newtheorem{thm}{Theorem}

\numberwithin{equation}{section}
\numberwithin{lem}{section}
\numberwithin{alg}{section}

\newcommand{\real}{\mathbb{R}}
\newcommand{\integer}{\mathbb{Z}}

\newcommand{\ddd}{\mathbb{D}}

\newcommand{\ttt}{\mathbb{T}}

\newcommand{\mlap}{\widetilde{\Delta}}

\DeclareMathOperator{\spn}{span}

\title{\vspace*{-1in}Orthogonal polynomials on the disk in the absence of finite moments}
\author{Peter C.~Gibson\footnote{Dept.~of Mathematics \& Statistics, York University, 4700 Keele St., Toronto, Ontario, Canada, M3J~1P3, $\mathtt{pcgibson@yorku.ca}$} }
\date{December 24, 2014}
\makeatletter
\let\newtitle\@title
\let\newauthor\@author
\let\newdate\@date
\makeatother

\begin{document}
\maketitle
\begin{abstract}
We introduce a new family of orthogonal polynomials on the disk that has emerged in the context of wave propagation in layered media.   Unlike known examples, the polynomials are orthogonal with respect to a measure all of whose even moments are infinite.   
\end{abstract}


\section{Introduction}

For each $\alpha>-1$ there is a corresponding family of disk polynomials that are orthogonal with respect to the measure $(1-x^2-y^2)^\alpha dxdy\,$ on the unit disk $\ddd$; these are sometimes referred to as generalized Zernike polynomials, named for the case  $\alpha=0$ introduced in \cite{Ze:1934}.   The well-established theory of disk polynomials is detailed in \cite{Ko:1975, Wu:2005, DuXu:2001}.  The constraint $\alpha>-1$ stems from the requirement that the measure $(1-x^2-y^2)^\alpha dxdy\,$ have finite moments, which is necessary for meaningful evaluation of the corresponding scalar product 
\begin{equation}\label{scalar}
\langle p,q\rangle_\alpha=\int_{\ddd}p(x,y)\overline{q(x,y)}(1-x^2-y^2)^\alpha dxdy
\end{equation}
on arbitrary polynomials $p$ and $q$.  Recent work on the propagation of waves in layered media \cite{Gi:SIAP2014, Gi:Smooth2014} has brought to light a family of polynomials orthogonal with respect to $(1-x^2-y^2)^{-1}dxdy$.  Since  
\begin{equation}\label{moments}
\int_{\ddd}x^{2m}y^{2n}(1-x^2-y^2)^{-1} dxdy=\infty
\end{equation}
for every pair of nonnegative integers $m$ and $n$, the scalar product (\ref{scalar}) is not defined for arbitrary polynomials in the case $\alpha=-1$.  Nevertheless, polynomials---which we term scattering polynomials---comprise an orthogonal basis for $L^2\bigl(\ddd,dxdy/(1-x^2-y^2)\bigr)$.  The purpose of the present paper is to present the details of this result.  

\section{Definition and properties of scattering polynomials\label{sec-main}}

Referring to the notation $z=x+iy$ for points in the unit disk $\ddd$, one has the option of working with euclidean $x,y$-coordinates, or with complex coordinates $z$ and $\bar{z}$.  As far as orthogonal polynomials are concerned these are essentially equivalent, as elaborated in \cite{Xu:2015}; the present paper uses whichever coordinates are most convenient for the task at hand.     

We define \emph{scattering polynomials} by a Rodrigues type formula, as follows.  For every $(p,q)\in\integer^2$ with min$\{p,q\}\geq 1$, set 
\begin{equation}\label{defn-scattering}
\varphi^{(p,q)}(z)=\frac{(-1)^p}{q(p+q-1)!}(1-z\bar{z})\frac{\partial^{p+q}}{\partial z^p\partial \bar{z}^q}(1-z\bar{z})^{p+q-1}.  
\end{equation}
The chosen normalization simplifies the formulation of the boundary Green's function for scattering in layered media (see \cite{Gi:Smooth2014}) and so is physically natural, although not important for present considerations.  

Note that disk polynomials satisfy a Rodrigues formula similar to that of scattering polynomials, but there is a qualitative difference: it follows directly from (\ref{defn-scattering}) that $\varphi^{(p,q)}(z)=0$ for every $z$ on the unit circle $\ttt$, whereas all disk polynomials have constant non-zero modulus on $\ttt$, cf.~\cite{Ko:1975}.   Our main result concerns completeness of scattering polynomials, as follows.  
\begin{thm}\label{thm-main}
Scattering polynomials $\varphi^{(p,q)}$ defined by (\ref{defn-scattering}), where $(p,q)\in\integer^2$ and $\min\{p,q\}\geq 1$, comprise an orthogonal basis for $L^2\bigl(\ddd,dxdy/(1-x^2-y^2)\bigr)$.  
\end{thm}
In \S\ref{sec-eigenfunctions} and \S\ref{sec-Jacobi} below we show that scattering polynomials are eigenfunctions of a second order differential operator and may be expressed in terms of Jacobi polynomials; these results contribute to a proof of Theorem~\ref{thm-main} completed in \S\ref{sec-completeness}.  

\subsection{Eigenfunctions of $-(1-x^2-y^2)\Delta/4$\label{sec-eigenfunctions}}

\pagestyle{fancyplain}

Let $\mlap$ denote the modified laplacian
\begin{equation}\label{mlap}
\mlap=(1-z\bar{z})\frac{\partial^2}{\partial z\partial \bar{z}}=\frac{1-x^2-y^2}{4}\Delta,
\end{equation}
where $\Delta$ is the usual (euclidean) laplacian.  Direct computation using (\ref{defn-scattering}) shows that for all integers $p,q\geq 1$, 
\begin{equation}\label{eigenfunctions}
-\mlap\varphi^{(p,q)}=pq\,\varphi^{(p,q)}.
\end{equation}
Letting $\sigma_0:\integer_+\rightarrow\integer_+$ denote the divisor function, there is thus a family of $\sigma_0(k)$ eigenfunctions of $-\mlap$ of the form $\varphi^{(p,q)}$ corresponding to each positive integer eigenvalue $k$.  We show in the next section that these eigenfunctions are linearly independent.  

\subsection{Representation in terms of {J}acobi polynomials\label{sec-Jacobi}}

Like disk polynomials, scattering polynomials have a representation in terms of Jacobi polynomials, but again, there is a qualitative difference.   The disk polynomials corresponding to parameter $\alpha>-1$ can be expressed in terms of Jacobi polynomials $P^{(\alpha,\beta)}_n$ for nonnegative integer values of $\beta$.   Since there is no Jacobi polynomial corresponding to $\alpha=-1$, the same cannot be true for scattering polynomials.  Indeed it turns out that scattering polynomials can be formulated in terms of $P^{(1,\beta)}_n$, where $\beta$ is a nonnegative integer, as follows.  

Expanding the binomial $(1-z\bar{z})^{p+q-1}$ in the formula (\ref{defn-scattering}), and then applying the derivative $\partial^{p+q}/\partial z^p\partial\bar{z}^q$, yields
\begin{equation}\label{step-one}
\varphi^{(p,q)}(z)=
\frac{(-1)^{q+\nu+1}}{q}(1-z\bar{z})z^{m+\nu-p+1}\bar{z}^{m+\nu-q+1}\sum_{j=0}^{\nu}(-1)^j\frac{(j+\nu+m+1)!}{j!(j+m)!(\nu-j)!}(z\bar{z})^j,
\end{equation}
where $m=|p-q|$ and $\nu=\min\{p,q\}-1$; the latter notation will be used in the remainder of this section.  Switching to polar form $z=re^{i\theta}$, it follows from (\ref{step-one}) that 
\begin{equation}\label{step-two}
\varphi^{(p,q)}\bigl(re^{i\theta}\bigr)=e^{i(q-p)\theta}f^{(p,q)}(r),
\end{equation}
where 
\begin{equation}\label{fpq}
f^{(p,q)}(r)=\frac{(-1)^{q+\nu+1}}{q}(1-r^2)r^m\sum_{j=0}^{\nu}(-1)^j\frac{(j+\nu+m+1)!}{j!(j+m)!(\nu-j)!}r^{2j}.
\end{equation}
The radial functions $f^{(p,q)}$ were first discovered in \cite{Gi:SIAP2014}, as was the following connection to Jacobi polynomials, valid for $\nu\geq 0$:   
\begin{equation}\label{Jacobi}
f^{(p,q)}(r)=\frac{(-1)^{q+m+\nu+1}(m+\nu+1)}{q}(1-r^2)r^mP^{(1,m)}_{\nu}(2r^2-1).  
\end{equation}
Combined with (\ref{step-two}) this yields the representation
\begin{equation}\label{representation}
\varphi^{(p,q)}\bigl(re^{i\theta}\bigr)=\frac{(-1)^{q+\max\{p,q\}}\max\{p,q\}}{q}(1-r^2)r^{|p-q|}P^{(1,|p-q|)}_{\min\{p,q\}-1}(2r^2-1)e^{i(q-p)\theta}.
\end{equation}

Note that the angular part of $\varphi^{(p,q)}(re^{i\theta})$, namely $e^{i(q-p)\theta}$, is a pure frequency.  Therefore if 
$
q-p\neq q^\prime-p^\prime,
$
then $\varphi^{(p,q)}$ and $\varphi^{(p^\prime,q^\prime)}$ are orthogonal in 
\[
L^2\Bigl(\ddd,\frac{dxdy}{1-x^2-y^2}\Bigr)=L^2\Bigl(\ddd,\frac{rdrd\theta}{1-r^2}\Bigr).
\]
In particular, if $pq=p^\prime q^\prime$ and $(p,q)\neq(p^\prime,q^\prime)$, then $\varphi^{(p,q)}$ and $\varphi^{(p^\prime,q^\prime)}$ are orthogonal, so the set of scattering polynomials corresponding to any fixed eigenvalue of $-\mlap$ is linearly independent.

\subsection{Completeness in $L^2\!\left(\ddd,\frac{r\,drd\theta}{1-r^2}\right)$\label{sec-completeness}}

In general, given a measure $\mu$ on a locally compact metric space $X$ and a positive measurable weight function $w:X\rightarrow\real_+$,
\begin{equation}\label{L2}
L^2(X,w\,d\mu)=\frac{1}{\sqrt{w}}L^2(X,d\mu),
\end{equation}
and a sequence $\{b_\nu\}_{\nu=0}^\infty$ is an orthogonal basis for $L^2(X,w\,d\mu)$ if and only if the corresponding sequence $\{\sqrt{w}b_\nu\}_{\nu=0}^\infty$ is an orthogonal basis for $L^2(X,d\mu)$.  
In particular,
setting $d\mu=rdrd\theta/(1-r^2)$, 
\begin{equation}\label{L2-relation}
L^2(\ddd,d\mu)=\sqrt{1-r^2}L^2(\ddd,rdrd\theta).
\end{equation}
Also, since for any nonnegative integer $m$, $\left\{P^{(1,m)}_\nu(u)\right\}_{\nu=0}^\infty$ is an orthogonal basis for
\[
L^2\bigl([-1,1],(1-u)(1+u)^m\,du\bigr),
\]
it follows that the quasipolynomials 
\begin{equation}\label{Q}
Q^{(1,m)}_\nu(u)=\left(\frac{1-u}{2}\right)^{\frac{1}{2}}\left(\frac{1+u}{2}\right)^{\frac{m}{2}}P^{(1,m)}_\nu(u)
\end{equation}
comprise an orthogonal basis for $L^2([-1,1],du)$; see \cite{Sz:1975}.   

In order to show that 
\begin{equation}\label{basis}
\mathcal{B}=\left\{\varphi^{(p,q)}\,\left|\,(p,q)\in\integer^2\;\&\;\min\{p,q\}\geq 1\right.\right\}
\end{equation}
is an orthogonal basis of $L^2\bigl(\ddd,rdrd\theta/(1-r^2)\bigr)$, we first argue that the functions $\varphi^{(p,q)}$ are orthogonal, and then that the span of $\mathcal{B}$ is dense.   It was proven in \S\ref{sec-Jacobi} that $\varphi^{(p,q)}$ and $\varphi^{(p^\prime,q^\prime)}$ are orthogonal if $pq=p^\prime q^\prime$ and $(p,q)\neq(p^\prime,q^\prime)$.  On the other hand, if $pq\neq p^\prime q^\prime$, then orthogonality of $\varphi^{(p,q)}$ and $\varphi^{(p^\prime,q^\prime)}$ follows from the fact that they are eigenfunctions, corresponding to distinct eigenvalues, of the self-adjoint operator $-\mlap$; self-adjointness of $-\mlap$ follows from that of $-\Delta$ by (\ref{mlap}).  

It remains to show that $\spn\mathcal{B}$ is dense in $L^2\bigl(\ddd,rdrd\theta/(1-r^2)\bigr)$.  Toward this end, suppose that $h\in L^2\bigl(\ddd,rdrd\theta/(1-r^2)\bigr)$ is orthogonal to every member of $\mathcal{B}$.  By (\ref{L2-relation}) there exists $g\in L^2\bigl(\ddd,rdrd\theta\bigr)$ such that 
\begin{equation}\label{gh}
h(r,\theta)=\sqrt{1-r^2}\,g(r,\theta).  
\end{equation}
Let $\alpha_{p,q}=(-1)^{q+\max\{p,q\}}\max\{p,q\}/q$ denote the coefficient occurring on the right-hand side of (\ref{representation}). Then for each fixed $n=q-p\in\integer$, for every $\nu=\min\{p,q\}-1\geq 0$,
\begin{equation}\nonumber
\begin{split}
0&=\int_\ddd h(r,\theta)\,\overline{\varphi^{(p,q)}(r,\theta)}\,\frac{rdrd\theta}{1-r^2}\\
 &=\alpha_{p,q}\int_0^1\left(\int_0^{2\pi}h(r,\theta)e^{-in\theta}\,d\theta\right)(1-r^2)r^{|n|}P^{(1,|n|)}_{\nu}(2r^2-1)\,\frac{rdr}{1-r^2}\qquad(\mbox{ by (\ref{Jacobi})})\\
&=\alpha_{p,q}\int_0^1\left(\int_0^{2\pi}g(r,\theta)e^{-in\theta}\,d\theta\right)\sqrt{1-r^2}\,r^{|n|}P^{(1,|n|)}_{\nu}(2r^2-1)\,rdr\qquad(\mbox{by (\ref{gh})})\\
&=\frac{\alpha_{p,q}}{4}\int_{-1}^1\left(\int_0^{2\pi}g\left(\textstyle\sqrt{\frac{1+u}{2}},\theta\right)e^{-in\theta}\,d\theta\right)\sqrt{\frac{1-u}{2}}\,\sqrt{\frac{1+u}{2}}^{|n|}P^{(1,|n|)}_{\nu}(u)\,du\qquad(u=2r^2-1)\\
&=\frac{\alpha_{p,q}}{4}\int_{-1}^1\left(\int_0^{2\pi}g\left(\textstyle\sqrt{\frac{1+u}{2}},\theta\right)e^{-in\theta}\,d\theta\right)Q^{(1,|n|)}_{\nu}(u)\,du\qquad(\mbox{ as in (\ref{Q})}).\\
\end{split}
\end{equation}
Since the quasipolynomials $Q^{(1,|n|)}_\nu$ are an orthogonal basis for $L^2([-1,1],du)$, it follows that for each $n\in\integer$,
\begin{equation}\label{g0}
\int_0^{2\pi}g\left(\textstyle\sqrt{\frac{1+u}{2}},\theta\right)e^{-in\theta}\,d\theta=0
\end{equation}
for every $u\in[-1,1]$ outside a set $E_n$ of measure zero.  Since $\{e^{in\theta}\}_{n\in\integer}$ is an orthogonal basis of $L^2([0,2\pi],d\theta)$, it follows in turn that for $u\not\in\cup E_n$
\[
 g\left(\textstyle\sqrt{\frac{1+u}{2}},\theta\right)=0,
 \]
for almost every $\theta\in[0,2\pi]$.  Thus $g(r,\theta)=0$ for almost every $(r,\theta)\in\ddd$ and $g=0$ as a function in $L^2(\ddd,rdrd\theta)$, whence $h=0$ also.  This proves that the orthogonal complement of $\mathcal{B}$ in $L^2\bigl(\ddd,rdrd\theta/(1-r^2)\bigr)$ is empty, and hence that $\mathcal{B}$ is an orthogonal basis.

\section{Conclusions}

Since the vector space $L^2\bigl(\ddd,rdrd\theta/(1-r^2)\bigr)=\sqrt{1-r^2}L^2\bigl(\ddd,rdrd\theta\bigr)$ is dense in 
\[
L^2\bigl(\ddd,rdrd\theta\bigr)=L^2\bigl(\ddd,dxdy\bigr),
\]
and convergence in the former space implies convergence in the latter, scattering polynomials comprise a (non-orthogonal) basis for$L^2\bigl(\ddd,dxdy\bigr)$ consistent with Dirichlet boundary values.  From the perspective of analysis of functions on the disk, this provides an alternative to Zernike polynomials---and their generalizations the disk polynomials---which are non-zero on the boundary circle and so inconsistent with Dirichlet conditions. 

More generally, scattering polynomials illustrate that orthogonal polynomials can comprise an orthogonal basis for a function spaces $L^2(X,d\mu)$ in which not all polynomials are integrable.  A natural question for further investigation is the existence and extent of other such examples.

\end{document}